\documentclass[11pt]{article}
\usepackage[]{amsmath,amssymb}

\newtheorem{theorem}{Theorem}[section]
\newtheorem{lemma}[theorem]{Lemma}

\newtheorem{exAux}[theorem]{Example}
\newenvironment{example}{\begin{exAux} \rm}{\end{exAux}}
\newtheorem{Def}[theorem]{Definition}
\newenvironment{definition}{\begin{Def} \rm}{\end{Def}}
\newtheorem{Rem}[theorem]{Remark}
\newenvironment{remark}{\begin{Rem} \rm}{\end{Rem}}
\newtheorem{Ass}[theorem]{Assumption}

\newenvironment{proof}{\medskip\noindent{\bf Proof.\ }}{\qed\medskip}
\newenvironment{proofof}[1]{\medskip\noindent{\bf Proof  of {#1}.\ 
}}{\qed\medskip}
\newcommand{\qed}{\hfill\mbox{$\Box$\qquad\qquad}}

\begin{document}
\thispagestyle{empty}
\vspace{1cm}

\begin{center}
\LARGE \bf
\noindent
Balanced Leonard Pairs
\end{center}

\begin{center}
\Large
Kazumasa Nomura and Paul Terwilliger
\end{center}

\smallskip

{\small 
\begin{center} {\bf Abstract}
\end{center}
Let $\mathbb{K}$ denote a field, and let $V$ denote a vector
space over $\mathbb{K}$ with finite positive dimension.
By a {\em Leonard pair} on $V$ we mean an ordered pair of linear transformations
$A:V \to V$ and $A^* : V \to V$
that satisfy the following two conditions:
\begin{itemize}
\item[(i)] There exists a basis for $V$ with respect to which
the matrix representing $A$ is irreducible tridiagonal and the
matrix representing $A^*$ is diagonal.
\item[(ii)] There exists a basis for $V$ with respect to which
the matrix representing $A^*$ is irreducible tridiagonal and the
matrix representing $A$ is diagonal.
\end{itemize}
Let $v^*_0$, $v^*_1$, \ldots, $v^*_d$ (respectively
$v_0$, $v_1$, \ldots, $v_d$)
denote a basis for $V$ that satisfies (i) (respectively (ii)). 
For $0 \leq i \leq d$, let  $a_i$ denote the coefficient of 
$v^*_i$, when we write
$A v^*_i$ as a linear combination of $v^*_0$, $v^*_1$, \ldots,
$v^*_d$, and let $a^*_i$ denote the coefficient of $v_i$, when we write
$A^* v_i$ as a linear combination of $v_0$, $v_1$, \ldots, $v_d$.

In this paper we show $a_0=a_d$ if and only if $a^*_0=a^*_d$.
Moreover we show that for $d \geq 1$ the following are equivalent;
(i) $a_0=a_d$ and $a_1=a_{d-1}$; (ii) $a^*_0=a^*_d$ and 
$a^*_1=a^*_{d-1}$;
(iii) $a_i=a_{d-i}$ and $a^*_i=a^*_{d-i}$ for $0 \leq i \leq d$.
These give a proof of a conjecture by the second author. 
We say $A$, $A^*$ is {\em balanced} whenever 
$a_i=a_{d-i}$ and $a^*_i=a^*_{d-i}$ for $0 \leq i \leq d$.
We say $A$, $A^*$ is {\em essentially bipartite} 
(respectively {\em essentially dual bipartite})
whenever $a_i$ (respectively $a^*_i$) is independent of $i$ for $0 \leq i \leq d$.
Observe that if $A$, $A^*$ is essentially bipartite or dual bipartite,
then $A$, $A^*$ is balanced. For $d \neq 2$, we show that if $A$, $A^*$ is
balanced then $A$, $A^*$ is essentially bipartite or dual bipartite.
}

\section{Introduction}

Let $\mathbb{K}$ denote a field, and let $V$ denote a vector
space over $\mathbb{K}$ with finite positive dimension.
We consider an ordered pair of linear transformations
$A:V \to V$ and $A^* : V \to V$
that satisfy the following two conditions:
\begin{itemize}
\item[(i)] There exists a basis for $V$ with respect to which
the matrix representing $A$ is irreducible tridiagonal and the
matrix representing $A^*$ is diagonal.
\item[(ii)] There exists a basis for $V$ with respect to which
the matrix representing $A^*$ is irreducible tridiagonal and the
matrix representing $A$ is diagonal.
\end{itemize}
Such a pair is called a {\em Leonard pair} on $V$. This notion was
introduced by the second author \cite{T:Leonard}.

\medskip
Throughout this paper, we fix the following notation.
Let $A$, $A^*$ denote a Leonard pair on $V$. We set $d= \dim V - 1$.
Let $v^*_0$, $v^*_1$, \ldots, $v^*_d$ denote a basis for $V$ that
satisfies the condition (i), and let
$v_0$, $v_1$, \ldots, $v_d$  denote a basis for $V$ that satisfies (ii).
For $0 \leq i \leq d$, let  $a_i$ denote the coefficient of 
$v^*_i$, when we write
$A v^*_i$ as a linear combination of $v^*_0$, $v^*_1$, \ldots,
$v^*_d$, and let $a^*_i$ denote the coefficient of $v_i$, when we write
$A^* v_i$ as a linear combination of $v_0$, $v_1$, \ldots, $v_d$.

\medskip

In this paper we prove the following results.

\medskip

\begin{theorem}   \label{thm:main0}
\samepage
The following are equivalent.
\begin{itemize}
\item[(i)] $a_0=a_d$,
\item[(ii)] $a^*_0=a^*_d$.
\end{itemize}
\end{theorem}   

\medskip

\begin{theorem}    \label{thm:main}
For $d \geq 1$ the following are equivalent.
\begin{itemize}
\item[(i)]  $a_0=a_d$ and $a_1=a_{d-1}$,
\item[(ii)] $a^*_0 = a^*_d$ and $a^*_1=a^*_{d-1}$,
\item[(iii)] $a_i=a_{d-i}$ and $a^*_i=a^*_{d-i}$ for $0 \leq i \leq d$.
\end{itemize}
\end{theorem}

\medskip

We say that $A$, $A^*$ is {\em balanced} whenever
$a_i=a_{d-i}$ and $a^*_i=a^*_{d-i}$ for $0 \leq i \leq d$.

\medskip

\begin{remark}
Theorems \ref{thm:main0} and \ref{thm:main} give a proof of a
conjecture by the second author \cite[Section 36]{T:survey}.
\end{remark}

\begin{remark}
Pascasio \cite[Corollary 4.3]{Pascasio} proved
Theorem \ref{thm:main0} for the Leonard pairs that come from a
$Q$-polynomial distance-regular graph.
\end{remark}

\medskip

For $0 \leq i \leq d$, let $\theta_i$ (respectively $\theta^*_i$)
denote the eigenvalue for $A$ associated with the eigenvector
$v_i$ (respectively $v^*_i$). 
Let $\varphi_1$, $\varphi_2$, \ldots, $\varphi_d$ 
(respectively $\phi_1$, $\phi_2$, \ldots, $\phi_d$) denote the
first split sequence (respectively the second split sequence)
with respect to the ordering 
$(\theta_0$, $\theta_1$, \ldots, $\theta_d$;
 $\theta^*_0$, $\theta^*_1$, \ldots, $\theta^*_d)$. The definition
of the split sequences will be given in Section 2.

\medskip
A Leonard pair is said to be {\em bipartite} whenever $a_i=0$
for $0 \leq i \leq d$. We consider a slightly more general situation.

\medskip

\begin{theorem}       \label{thm:X}
The following are equivalent.
\begin{itemize}
\item[(i)] $a_i$ is independent of $i$ for $0 \leq i \leq d$.
\item[(ii)] $\theta_i + \theta_{d-i}$ is independent of $i$
for $0 \leq i \leq d$, and $\varphi_i = - \phi_i$ for $1 \leq i \leq d$.
\end{itemize}
Suppose (i), (ii) hold. Then the common value of
$\theta_i + \theta_{d-i}$ is twice the common value of the $a_i$.
\end{theorem}

\medskip
We say the Leonard pair $A$, $A^*$ is {\em essentially bipartite} whenever the
equivalent conditions (i), (ii) hold in Theorem \ref{thm:X}.
Observe that if  $A$, $A^*$ is essentially bipartite, then
the Leonard pair $A - \xi I$, $A^*$ is bipartite, where
$\xi$ denotes the common value of $a_0$, $a_1$, \ldots, $a_d$.

\medskip
A Leonard pair is said to be {\em dual bipartite} whenever $a^*_i=0$
for $0 \leq i \leq d$. We consider a slightly more general situation.

\medskip

\begin{theorem}       \label{thm:Xs}
\samepage
The following are equivalent.
\begin{itemize}
\item[(i)] $a^*_i$ is independent of $i$ for $0 \leq i \leq d$.
\item[(ii)] $\theta^*_i + \theta^*_{d-i}$ is independent of $i$
for $0 \leq i \leq d$, and $\varphi_i = - \phi_{d-i+1}$ for $1 \leq i \leq d$.
\end{itemize}
Suppose (i), (ii) hold. Then the common value of
$\theta^*_i + \theta^*_{d-i}$ is twice the common value of the $a^*_i$.
\end{theorem}

\medskip
We say the Leonard pair $A$, $A^*$ is {\em essentially dual bipartite} whenever the
equivalent conditions (i), (ii) hold in Theorem \ref{thm:Xs}.
Observe that if $A$, $A^*$ is essentially dual bipartite, then
the Leonard pair $A$, $A^* - \xi^* I$ is dual bipartite, where
$\xi^*$ denotes the common value of $a^*_0$, $a^*_1$, \ldots, $a^*_d$.

\medskip

\begin{theorem}        \label{thm:Y}
Let $A$, $A^*$ denote a Leonard pair.
\begin{itemize}
\item[(i)] If $A$, $A^*$ is essentially bipartite, then $A$, $A^*$ is balanced.
\item[(ii)] If $A$, $A^*$ is essentially dual bipartite, then 
$A$, $A^*$ is balanced.
\item[(iii)] Assume $d \neq 2$. 
If $A$, $A^*$ is balanced, then $A$, $A^*$ is essentially bipartite or 
essentially dual bipartite.
\end{itemize}
\end{theorem}

\medskip

\begin{remark}
For $d=2$, part (iii) of Theorem \ref{thm:Y} is false. 
A counter example is given in Example \ref{exm:counter}.
\end{remark}

\medskip

\begin{remark}
In our proof of Theorems \ref{thm:main}, \ref{thm:X}--\ref{thm:Y}
we use a case-analysis based on the classification
of Leonard pairs by the second author \cite{T:Leonard,T:array}.
\end{remark}

\medskip

Our paper is organized as follows.
In Section 2 we give some background information.
In Section 3 we give the proof of Theorem \ref{thm:main0}.
In Section 4 we describe the cases that we will use in our proof
of Theorems \ref{thm:main}, \ref{thm:X}--\ref{thm:Y}.
In Sections 5--10 we give the proofs of these theorems.

\section{Some background information}

In this section we summarize some results that we will use
in our proof.

\medskip

\begin{lemma} \cite[Lemma 1.3]{T:Leonard}   \label{lem:eigen}
The eigenvalues $\theta_0$, $\theta_1$, \ldots, $\theta_d$ of $A$
are distinct and contained in $\mathbb{K}$. Moreover, the eigenvalues
$\theta^*_0$, $\theta^*_1$, \ldots, $\theta^*_d$ of $A^*$ are distinct and
contained in $\mathbb{K}$.
\end{lemma}

\medskip

\begin{lemma}   \cite[Lemma 9.5]{T:Leonard}     \label{lem:0th}
For $d \geq 1$ and for $0 \leq i \leq d$,
\begin{equation}         \label{eq:main0aux}
    \frac{\theta_i-\theta_{d-i}}{\theta_0 - \theta_d}
  = \frac{\theta^*_i -\theta^*_{d-i}}{\theta^*_0-\theta^*_{d}}.
\end{equation}
\end{lemma}

\medskip

\begin{theorem} \cite[Theorem 3.2]{T:Leonard}   \label{thm:split}
There exists a basis for $V$ with
respect to which the matrices representing $A$, $A^*$ take the
following form for some scalars $\varphi_1$, $\varphi_2$, \ldots, $\varphi_d$ in
$\mathbb{K}$:
\[
A: 
\begin{pmatrix}
   \theta_0 &          & & & & \text{\bf 0} \\
   1        & \theta_1 \\
            & 1        & \theta_2 \\
            &          &  \cdot   &  \cdot \\
            &          &          & \cdot  & \cdot\\
   \text{\bf 0} &      &          &       & 1 & \theta_d
\end{pmatrix},
\qquad
A^*:
\begin{pmatrix}
   \theta^*_0 &  \varphi_1  & & & & \text{\bf 0} \\
           & \theta^*_1 & \varphi_2 \\
            &         & \theta^*_2 & \cdot \\
            &          &       &  \cdot & \cdot \\
            &          &          &   & \cdot & \varphi_d\\
   \text{\bf 0} &      &          &       &  & \theta^*_d
\end{pmatrix}.
\]
The sequence $\varphi_1$, $\varphi_2$, \ldots, $\varphi_d$ is
uniquely determined by the ordering 
$(\theta_0$, $\theta_1$, \ldots, $\theta_d$;
$\theta^*_0$, $\theta^*_1$, \ldots, $\theta^*_d)$.
Moreover $\varphi_i \neq 0$ for $1 \leq i \leq d$.
\end{theorem}

\medskip

The sequence $\varphi_1$, $\varphi_2$, \ldots, $\varphi_d$ is called the
{\em first split sequence} with respect to the ordering
$(\theta_0$, $\theta_1$, \ldots, $\theta_d$; 
$\theta^*_0$, $\theta^*_1$, \ldots, $\theta^*_d)$.
Let $\phi_1$, $\phi_2$, \ldots, $\phi_d$ denote the first
split sequence with respect to the ordering
$(\theta_d$, $\theta_{d-1}$, \ldots, $\theta_0$; 
$\theta^*_0$, $\theta^*_1$, \ldots, $\theta^*_d)$.
We call $\phi_1$, $\phi_2$, \ldots, $\phi_d$ the
{\em second split sequence} with respect to the ordering
$(\theta_0$, $\theta_1$, \ldots, $\theta_d$; 
$\theta^*_0$, $\theta^*_1$, \ldots, $\theta^*_d)$.
The sequence
\[
   (\theta_0, \theta_1, \ldots, \theta_d; \;
   \theta^*_0, \theta^*_1, \ldots, \theta^*_d; \;
   \varphi_1, \varphi_2, \ldots, \varphi_d; \;
   \phi_1, \phi_2, \ldots, \phi_d)
\]
is called a {\em parameter array} of the Leonard pair.

\medskip
In the classification of Leonard pairs, the following
theorem plays a key role.

\medskip

\begin{theorem}   \cite[Theorem 1.9]{T:Leonard}     \label{thm:classify}
Let 
\begin{equation}              \label{eq:paramarray}
   (\theta_0, \theta_1, \ldots, \theta_d; \;
   \theta^*_0, \theta^*_1, \ldots, \theta^*_d; \;
   \varphi_1, \varphi_2, \ldots, \varphi_d; \;
   \phi_1, \phi_2, \ldots, \phi_d)
\end{equation}
denote a sequence of scalars taken from $\mathbb{K}$. 
Then there exists a Leonard pair with
parameter array (\ref{eq:paramarray}) if and only if (i)--(v) hold
below.
\begin{itemize}
\item[(i)]  $\varphi_i \neq 0$, $\phi_i \neq 0$ $(1 \leq i \leq d)$.
\item[(ii)] $\theta_i \neq \theta_j$, $\theta^*_i \neq \theta^*_j$
   if $i \neq j$ $(0 \leq i,j \leq d$).
\item[(iii)] For $1 \leq i \leq d$,
\[
 \varphi_i = \phi_1 
 \sum_{h=0}^{i-1} \frac{\theta_h - \theta_{d-h}}
                       {\theta_0 - \theta_d}
         + (\theta^*_{i}-\theta^*_{0})(\theta_{i-1}-\theta_{d}).
\]
\item[(iv)] For $1 \leq i \leq d$,
\[
\phi_i = \varphi_1
\sum_{h=0}^{i-1} \frac{\theta_h - \theta_{d-h}}
                       {\theta_0 - \theta_d}
         + (\theta^*_{i}-\theta^*_{0})(\theta_{d-i+1}-\theta_{0}).
\]
\item[(v)] The expressions
\begin{equation}        \label{eq:indep}
   \frac{\theta_{i-2}-\theta_{i+1}}{\theta_{i-1}-\theta_{i}},
 \quad
   \frac{\theta^*_{i-2}-\theta^*_{i+1}}{\theta^*_{i-1}-\theta^*_{i}}
\end{equation}
are equal and independent of $\;i\;$ for $2 \leq i \leq d-1$.
\end{itemize}
\end{theorem}

\medskip

The scalars $a_i$, $a^*_i$ can be expressed in terms of the
parameter array as follows.

\medskip

\begin{lemma} \cite[Lemma 10.3]{T:24points} 
For $0 \leq i \leq d$, 
\begin{equation}         \label{eq:ai}
 a_i = \theta_i 
       + \frac{\varphi_i}{\theta^*_i-\theta^*_{i-1}}
       - \frac{\varphi_{i+1}}{\theta^*_{i+1}-\theta^*_i},  \qquad 
 a^*_i = \theta^*_i 
       + \frac{\varphi_i}{\theta_i - \theta_{i-1}}
       - \frac{\varphi_{i+1}}{\theta_{i+1} - \theta_i},
\end{equation}
\begin{equation}          \label{eq:ai2}     
 a_i = \theta_{d-i} 
       + \frac{\phi_i}{\theta^*_i-\theta^*_{i-1}}
       - \frac{\phi_{i+1}}{\theta^*_{i+1}-\theta^*_i},  \qquad 
 a^*_i = \theta^*_{d-i} 
       + \frac{\phi_{d-i+1}}{\theta_i - \theta_{i-1}}
       - \frac{\phi_{d-i}}{\theta_{i+1} - \theta_i},
\end{equation}
where we set $\varphi_0=0$, 
$\varphi_{d+1}=0$, $\phi_{0}=0$, $\phi_{d+1}=0$, and
let $\theta_{-1}$, $\theta_{d+1}$, $\theta^*_{-1}$,
$\theta^*_{d+1}$ denote indeterminates.
\end{lemma}

\section{Proof of Theorem \ref{thm:main0}}

In this section we prove Theorem \ref{thm:main0}.

\medskip

\begin{lemma}        \label{lem:a0andad}
For $d \geq 1$,
\begin{eqnarray}
a_0 &=& \theta_0 + \frac{\varphi_1}{\theta^*_0 - \theta^*_1},
       \label{eq:a0} \\
a_d &=& \frac{\theta_1(\theta^*_0-\theta^*_{d})-\theta_0(\theta^*_0-\theta^*_{d-1})}
           {\theta^*_{d-1}-\theta^*_d}
    - \frac{\varphi_1}{\theta^*_{d-1}-\theta^*_d},
       \label{eq:ad}  \\
a^*_0&=& \theta^*_0 + \frac{\varphi_1}{\theta_0 - \theta_1}
       \label{eq:as0}, \\
a^*_d &=&  \frac{\theta^*_1(\theta_0-\theta_{d})-\theta^*_0(\theta_0-\theta_{d-1})}
           {\theta_{d-1}-\theta_d}
             - \frac{\varphi_1}{\theta_{d-1}-\theta_d}.
       \label{eq:asd}  
\end{eqnarray}
\end{lemma}

\begin{proof}
The equations (\ref{eq:a0}) and (\ref{eq:as0}) follow from
(\ref{eq:ai}). 
From Theorem \ref{thm:classify} (iii), (iv),
\begin{equation}              \label{eq:varphid}
\varphi_d =  \varphi_1 + (\theta^*_1-\theta^*_0)(\theta_d-\theta_0)
             +  (\theta^*_d-\theta^*_0)(\theta_{d-1}-\theta_d).
\end{equation}
From (\ref{eq:main0aux}) at $i=1$,
\begin{equation}        \label{eq:thd}
  \theta_d = \theta_0
      - \frac{(\theta^*_0-\theta^*_d)(\theta_1-\theta_{d-1})}
             {\theta^*_1-\theta^*_{d-1}}.
\end{equation}
Evaluating the equation on the left in (\ref{eq:ai}) 
using (\ref{eq:varphid}) and (\ref{eq:thd})  we find
(\ref{eq:ad}).
The proof of (\ref{eq:asd}) is similar.
\end{proof}

\medskip

\begin{lemma}
For $d \geq 1$,
\begin{eqnarray}       
a_0-a_d &=& 
    \frac{(\theta_0-\theta_1)(\theta^*_0-\theta^*_d)}
         {\theta^*_{d-1} - \theta^*_d}    
  + \frac{\varphi_1}{\theta^*_0-\theta^*_1} 
  + \frac{\varphi_1}{\theta^*_{d-1}-\theta^*_d},
         \label{eq:a0ad}  \\
a^*_0-a^*_d &=&
    \frac{(\theta^*_0-\theta^*_1)(\theta_0-\theta_d)}
         {\theta_{d-1} - \theta_d}
  + \frac{\varphi_1}{\theta_0-\theta_1} 
  + \frac{\varphi_1}{\theta_{d-1}-\theta_d}.
       \label{eq:as0asd} 
\end{eqnarray}
\end{lemma}

\begin{proof}
Follows from Lemma \ref{lem:a0andad}.
\end{proof}

\medskip

\begin{lemma}    \label{lem:main0coincide}
For $d \geq 1$,
\begin{equation}       \label{eq:main0coincide}
  \frac{(a_0-a_d)(\theta^*_0-\theta^*_1)(\theta^*_{d-1}-\theta^*_d)}
       {\theta^*_0 - \theta^*_d}
= \frac{(a^*_0-a^*_d)(\theta_0-\theta_1)(\theta_{d-1}-\theta_d)}
       {\theta_0 - \theta_d}.
\end{equation}
\end{lemma}

\medskip

\begin{proof}
Using (\ref{eq:a0ad}) and (\ref{eq:as0asd}), 
the left side of (\ref{eq:main0coincide}) becomes
\[
  \varphi_1 + (\theta_0-\theta_1)(\theta^*_0-\theta^*_1)
           - \frac{\varphi_1(\theta^*_1-\theta^*_{d-1})}
                  {\theta^*_0 - \theta^*_d},
\]
and the right side of  (\ref{eq:main0coincide}) becomes
\[
  \varphi_1 + (\theta^*_0-\theta^*_1)(\theta_0-\theta_1)
           - \frac{\varphi_1(\theta_1-\theta_{d-1})}
                  {\theta_0 - \theta_d}.
\]
These expressions coincide by (\ref{eq:main0aux}).
\end{proof}

\medskip

\begin{proofof}{Theorem \ref{thm:main0}}
Assume $d \geq 1$; otherwise the result is vacuously true. Now the
result follows from Lemma \ref{lem:main0coincide}.
\end{proofof}

\section{Description of the cases}

Let $\overline{\mathbb{K}}$ denote the algebraic closure of
$\mathbb{K}$. 
In our proof of Theorems \ref{thm:main}, \ref{thm:X}, \ref{thm:Xs} and
\ref{thm:Y},
we break the argument into the following cases.

\medskip

\begin{itemize}
\item[] Case 0: $d \leq 2$.
\end{itemize}
For $d \geq 3$ let $q$ denote  a nonzero scalar in $\overline{\mathbb{K}}$
such that $q+q^{-1}+1$ is equal to the common value of (\ref{eq:indep}).
\begin{itemize}
\item[] Case I:  $d \geq 3$, $q \neq 1$, $q \neq -1$.
\item[] Case II: $d \geq 3$, $q=1$, $\text{\rm Char}(\mathbb{K}) \neq 2$.
\item[] Case III: $d \geq 3$, $q=-1$, $\text{\rm Char}(\mathbb{K}) \neq 2$, $d$ even.
\item[] Case IV: $d \geq 3$, $q=-1$, $\text{\rm Char}(\mathbb{K}) \neq 2$, $d$ odd.  
\item[] Case  V: $d \geq 3$, $q=1$, $\text{\rm Char}(\mathbb{K}) = 2$.
\end{itemize}

\medskip

\begin{definition}           \label{def:H}
For $d \geq 1$ we let $H$ denote the value of (\ref{eq:main0coincide});
\[
H=  \frac{(a_0-a_d)(\theta^*_0-\theta^*_1)(\theta^*_{d-1}-\theta^*_d)}
       {\theta^*_0 - \theta^*_d}
= \frac{(a^*_0-a^*_d)(\theta_0-\theta_1)(\theta_{d-1}-\theta_d)}
       {\theta_0 - \theta_d}.
\]
We note that $H=0$ if and only if $a_0=a_d$ if and only if
$a^*_0=a^*_d$.
\end{definition}

\section{Case 0: $d \leq 2$}

In this section we prove Theorems \ref{thm:main}, \ref{thm:X}--\ref{thm:Y} 
for $d \leq 2$.
We first note that Theorem \ref{thm:main} follows from Theorem \ref{thm:main0}
for these values of $d$.
We consider Theorems \ref{thm:X}--\ref{thm:Y}.

\medskip
First assume $d=0$. Then Theorems \ref{thm:X}, \ref{thm:Xs} and \ref{thm:Y} 
are vacuously true.

\medskip

Next assume $d=1$. 
From (\ref{eq:a0}) and the equation on the left in (\ref{eq:ai2}) for $i=1$,
\[
   a_0 - a_1 = \frac{\varphi_1 + \phi_1}{\theta^*_0 - \theta^*_1}.
\]
Thus $a_0=a_1$ if and only if $\varphi_1+\phi_1=0$.
From (\ref{eq:a0}) and (\ref{eq:ad}), we find $a_0+a_1 = \theta_0+\theta_1$.
So that $2 a_0 = \theta_0+\theta_1$ when $a_0=a_1$. 
These imply Theorem \ref{thm:X}.
The proof of Theorem  \ref{thm:Xs} is similar.
Theorem \ref{thm:Y} follows from Theorem \ref{thm:main0}.

\bigskip

For the rest of this section, we assume  $d=2$.

\medskip

\begin{lemma}      \label{lem:d2a0a2}
The following hold.
\begin{eqnarray}
 \varphi_1 &=& H -  (\theta_0-\theta_1)(\theta^*_0-\theta^*_1),
    \label{eq:d2vphi1}  \\
 \varphi_2 &=& H-(\theta_1-\theta_2)(\theta^*_1 - \theta^*_2),
    \label{eq:d2vphi2}  \\
 \phi_1 &=& H + (\theta_1-\theta_2)(\theta^*_0 - \theta^*_1),
    \label{eq:d2phi1}  \\
 \phi_2 &=& H + (\theta_0-\theta_1)(\theta^*_1-\theta^*_2).
    \label{eq:d2phi2}
\end{eqnarray}
\end{lemma}

\begin{proof}
Setting $d=2$ in (\ref{eq:a0ad}) we find (\ref{eq:d2vphi1}).
The other equations follow from (\ref{eq:d2vphi1}) using
Theorem \ref{thm:classify} (iii), (iv).
\end{proof}

\medskip

\begin{lemma}      \label{lem:d2a0a1}
Suppose $H=0$. Then
\[
      a_1-a_0 = \theta_0 - 2 \theta_1 + \theta_2.
\]
\end{lemma}

\begin{proof}
Obtained by evaluating the equation on the left in (\ref{eq:ai}) for $i=0$, $1$
using (\ref{eq:d2vphi1}) and (\ref{eq:d2vphi2}).
\end{proof}

\medskip

\begin{proofof}{Theorem \ref{thm:X}}
(i)$\Rightarrow$(ii):
By assumption $a_0=a_2$ so $H=0$.
Using $H=0$ and $a_0 = a_1$ we find $\theta_0 + \theta_2=2 \theta_1$ by
Lemma \ref{lem:d2a0a1}.
Evaluating the data in Lemma \ref{lem:d2a0a2} using these equations
we find $\varphi_1=-\phi_1$ and $\varphi_2 = -\phi_2$.

(ii)$\Rightarrow$(i):
Observe Char$(\mathbb{K}) \neq 2$; otherwise the equation 
$\theta_0+\theta_2=2\theta_1$ becomes $\theta_0=\theta_2$ for a
contradiction.
Comparing (\ref{eq:d2vphi1}), (\ref{eq:d2phi1}) we find $2H=0$ so $H=0$.
By this and Definition \ref{def:H} we find $a_0=a_2$.
Evaluating Lemma \ref{lem:d2a0a1} using $H=0$ and 
$\theta_0+\theta_2=2\theta_1$ we find $a_0=a_1$. Now $a_0=a_1=a_2$
as desired.

Suppose (i), (ii) hold. Evaluating (\ref{eq:a0}) using (\ref{eq:d2vphi1})
we find $a_0=\theta_1$, so that the common value of $\theta_i+\theta_{d-i}$
is $2 a_0$.
\end{proofof}

\medskip
\begin{proofof}{Theorem \ref{thm:Xs}}
Similar to the proof of Theorem \ref{thm:X}.
\end{proofof}

\medskip
\begin{proofof}{Theorem \ref{thm:Y}}
Follows from Theorem \ref{thm:main}.
\end{proofof}

\medskip

We finish this section by giving an example that shows Theorem \ref{thm:Y} (iii)
is false for $d=2$.

\medskip

\begin{example}         \label{exm:counter}
Let $\theta_0$, $\theta_1$, $\theta_2$, $\theta^*_0$, $\theta^*_1$, $\theta^*_2$
denote  scalars in $\mathbb{K}$ such that
$\theta_i \neq \theta_j$, $\theta^*_i \neq \theta^*_j$ if $i\neq j$
$(0 \leq i,j \leq 2)$.
We define scalars
\begin{eqnarray*}
 \varphi_1 &=& - (\theta_0-\theta_1)(\theta^*_0 - \theta^*_1),  \\
 \varphi_2 &=& - (\theta_1-\theta_2)(\theta^*_1 - \theta^*_2), \\
 \phi_1    &=&   (\theta_1-\theta_2)(\theta^*_0 - \theta^*_1),  \\
 \phi_2    &=&   (\theta_0-\theta_1)(\theta^*_1 - \theta^*_2).
\end{eqnarray*}
Observe that the sequence
\begin{equation}   \label{eq:counter}
   (\theta_0,\theta_1,\theta_2; \; \theta^*_0,\theta^*_1,\theta^*_2; \;
    \varphi_1,\varphi_2; \; \phi_1,\phi_2)
\end{equation}
satisfies the conditions (i)--(v) in Theorem \ref{thm:classify}, so that
there exists a Leonard pair having the parameter array (\ref{eq:counter}).
Using (\ref{eq:ai}), we get
\[
  a_0 = \theta_1,  \quad
  a_1 = \theta_0 - \theta_1 + \theta_2,  \quad
  a_2 = \theta_1,
\]
\[
  a^*_0 = \theta^*_1,  \quad
  a^*_1 = \theta^*_0 - \theta^*_1 + \theta^*_2,  \quad
  a^*_2 = \theta^*_1.
\]
Observe $a_0=a_2$ and $a^*_0=a^*_2$, so that the Leonard pair is
balanced. 
On the other hand, it is essentially bipartite if and only if
$\theta_1 = \theta_0 - \theta_1 + \theta_2$, and it is essentially
dual bipartite if and only if
$\theta^*_1 = \theta^*_0 - \theta^*_1 + \theta^*_2$.
Therefore it is not essentially bipartite, and is not
essentially dual biparitite for $2\theta_1 \neq \theta_0+\theta_2$ and
$2 \theta^*_1 \neq \theta^*_0 + \theta^*_2$.
\end{example}

\section{Case I: $d \geq 3$, $q\neq 1$, $q\neq -1$}

In this section we assume $d \geq 3$, $q \neq 1$, $q \neq -1$.

\medskip

\begin{theorem} \cite{T:array}        \label{lem:I}
There exist scalars $\eta$, $\mu$, $h$, $\eta^*$, $\mu^*$, $h^*$, $\tau$ in 
$\overline{\mathbb{K}}$ such that for $0 \leq i \leq d$
\begin{eqnarray}
\theta_i &=& \eta + \mu q^i + h q^{d-i},  
   \label{eq:Ith} \\
\theta^*_i &=& \eta^* + \mu^* q^i + h^* q^{d-i},
   \label{eq:Iths}
\end{eqnarray}
and for $1 \leq i \leq d$
\begin{eqnarray}
\varphi_i &=& (q^i-1)(q^{d-i+1}-1)(\tau -\mu \mu^* q^{i-1} - hh^* q^{d-i}),
   \label{eq:Iphi}  \\
\phi_i &=& (q^i-1)(q^{d-i+1}-1)(\tau -h \mu^* q^{i-1} - \mu h^* q^{d-i}).
   \label{eq:Iphi2}
\end{eqnarray}
\end{theorem}

\medskip

\begin{proof}
These are (27), (28), (31), (32) in \cite{T:array} 
after a change of variables.
\end{proof}

\medskip

\begin{remark}         \label{rem:I}
For $1 \leq i \leq d$ we have $q^i \neq 1$;
otherwise $\varphi_i=0$ by (\ref{eq:Iphi}).
\end{remark}

\smallskip

\begin{lemma}         \label{lem:Iada0}
$H = (q-1)^2((q^{d-1}+1) \tau - q^{d-1} (h+\mu)(h^* + \mu^*))$.
\end{lemma}

\begin{proof}
It is routine to verify this equation
using (\ref{eq:ai}), (\ref{eq:Ith}), (\ref{eq:Iths}) and (\ref{eq:Iphi}).
\end{proof}

\medskip

\begin{lemma}     \label{lem:Inonzero}
Assume $H=0$. Then $q^{d-1}+1 \neq 0$ and
\begin{equation}       \label{eq:Itau}
  \tau = \frac{q^{d-1} (h+\mu)(h^*+\mu^*)}{q^{d-1}+1}.
\end{equation}
\end{lemma}

\begin{proof}
Assume $q^{d-1}+1 = 0$. Then $1=-q^{d-1}$, so that
\begin{eqnarray*}
  0 &=& (q^{d-1}+1)\tau \\
    &=&  q^{d-1}(h + \mu)(h^* + \mu^*) \\
    &=&  q^{d-1}(\mu - h q^{d-1})(\mu^* - h^* q^{d-1}) \\
    &=&  q^{d-1} (q-1)^{-2}(\theta_0-\theta_1)(\theta^*_0-\theta^*_1),
\end{eqnarray*}
a contradiction, so we must have $q^{d-1}+1 \neq 0$
and (\ref{eq:Itau}) follows.
\end{proof}

\medskip

\begin{lemma}         \label{lem:Iadm1a1}
Assume $H=0$. Then the following coincide.
\[
    \frac{(a_{1}-a_{d-1})(\theta^*_0-\theta^*_3)(\theta^*_{d-3}-\theta^*_d)}
         {\theta^*_0 - \theta^*_d},
\]
\[
    \frac{(a^*_{1}-a^*_{d-1})(\theta_0-\theta_3)(\theta_{d-3}-\theta_d)}
         {\theta_0 - \theta_d},
\]
\[
  \frac{(1-q^2)(q^3-1)^2(q^{d-1}-1)(q^{d-2}-1)\tau}
       {q^2(q^d-1)}.
\]
\end{lemma}

\begin{proof}
It is routine to verify the coincidence
using (\ref{eq:ai}), (\ref{eq:Ith}), (\ref{eq:Iths}), (\ref{eq:Iphi})
and (\ref{eq:Itau}).
\end{proof}

\medskip

\begin{lemma}         \label{lem:Iequiv1}
Assume $H=0$. Then the following are equivalent.
\begin{itemize}
\item[(i)] $a_{1} = a_{d-1}$,
\item[(ii)] $a^*_1 = a^*_{d-1}$,
\item[(iii)] $(h+\mu)(h^* + \mu^*) = 0$,
\item[(iv)] $\tau = 0$.
\end{itemize}
\end{lemma}

\begin{proof}
Follows from Lemma \ref{lem:Iadm1a1} and  (\ref{eq:Itau}).
\end{proof}

\medskip

\begin{theorem}         \label{thm:I}
Assume $d \geq 3$, $q \neq 1$, $q \neq -1$.
Then the following are equivalent.
\begin{itemize}
\item[(i)] $a_0=a_d$ and $a_1=a_{d-1}$,
\item[(ii)] $a^*_0=a^*_d$ and $a^*_1 = a^*_{d-1}$,
\item[(iii)] $a_i = a_{d-i}$ and $a^*_i=a^*_{d-i}$ for $0 \leq i \leq d$,
\item[(iv)] $\tau=0$ and $(h+\mu)(h^*+\mu^*)=0$.
\end{itemize}
\end{theorem}

\begin{proof}
The conditions (i), (ii), (iv) are equivalent by Lemmas 
\ref{lem:Iada0} and \ref{lem:Iequiv1}. 
Clearly (iii) implies (i).
We show (iv) implies (iii).
Observe that we have $h= -\mu$ or $h^*=- \mu^*$.
For the case $h=-\mu$,
it is routine to verify $a_{d-i}-a_{i}=0$ and $a^*_{d-i}-a^*_{i}=0$
using (\ref{eq:ai}), (\ref{eq:Ith}), (\ref{eq:Iths}) and (\ref{eq:Iphi}) with
$\tau=0$ and $h=-\mu$. The case $h^*=-\mu^*$ is
similar.
\end{proof}

\medskip

\begin{lemma}          \label{lem:Igammagammas}
For $0 \leq i \leq d$,
\begin{eqnarray*}
  \theta_i+\theta_{d-i} &=& 2\eta + (h+\mu)(q^i+q^{d-i}), \\
  \theta^*_i+\theta^*_{d-i} &=& 2\eta^* + (h^*+\mu^*)(q^i+q^{d-i}).
\end{eqnarray*}
\end{lemma}

\begin{proof}
It is routine to verify these equations using
(\ref{eq:Ith}) and (\ref{eq:Iths}).
\end{proof}

\medskip

\begin{lemma}              \label{lem:Ivphiphi}
For $1 \leq i \leq d$,
\begin{eqnarray*}
  \varphi_i + \phi_i 
  &=& (q^i-1)(q^{d-i+1}-1)(2\tau - (h+\mu)(\mu^* q^{i-1}+h^* q^{d-i})), \\
  \varphi_i + \phi_{d-i+1}
 &=& (q^i-1)(q^{d-i+1}-1)(2\tau - (h^*+\mu^*)(\mu q^{i-1}+h q^{d-i})).
\end{eqnarray*}
\end{lemma}

\begin{proof}
It is routine to verify these equations using
(\ref{eq:Iphi}) and (\ref{eq:Iphi2}).
\end{proof}

\medskip

\begin{lemma}          \label{lem:Ia1a0}
The following hold.
\begin{itemize}
\item[(i)] Assume $\tau=0$ and $h^* + \mu^* = 0$.  Then
\[
a_1-a_0 =
   \frac{q^{d-2}(q-1)(q^2-1)^2(q^{d-1}-1)(\mu^*)^2(h+\mu)}
        {(\theta^*_0 - \theta^*_1)(\theta^*_1 - \theta^*_2)}.
\]
\item[(ii)] Assume $\tau=0$ and $h + \mu = 0$ then 
\[
a^*_1-a^*_0 =
   \frac{q^{d-2}(q-1)(q^2-1)^2(q^{d-1}-1)\mu^{2}(h^*+\mu^*)}
        {(\theta_0 - \theta_1)(\theta_1 - \theta_2)}.
\]
\end{itemize}
\end{lemma}

\begin{proof}
It is routine to verify these equations using
(\ref{eq:ai}), (\ref{eq:Ith}), (\ref{eq:Iths}), (\ref{eq:Iphi}),
(\ref{eq:Iphi2}).
\end{proof}

\begin{theorem}        \label{thm:Iaux}
Assume $d \geq 3$, $q \neq 1$, $q \neq -1$.
Then the following are equivalent.
\begin{itemize}
\item[(i)] $\tau=0$ and $h+\mu=0$.
\item[(ii)] $a_i$ is independent of $i$ for $0 \leq i \leq d$.
\item[(iii)] $\theta_i + \theta_{d-i}$ is independent of $i$
   for $0 \leq i \leq d$, and
   $\varphi_i = -\phi_i$ for $1 \leq i \leq d$.
\end{itemize}
Suppose (i)--(iii) hold. Then the common value of 
$a_i$ is $\eta$, and the common value of
$\theta_i + \theta_{d-i}$ is $2 \eta$.
\end{theorem}

\begin{proof}
(i)$\Rightarrow$(ii): Evaluating (\ref{eq:Ith}), (\ref{eq:Iphi}) 
using $\tau=0$ and $h=-\mu$ we find
\begin{eqnarray}
 \theta_i &=& \eta + \mu(q^i - q^{d-i}) \qquad (0 \leq i \leq d),   
     \label{eq:Iaux1}  \\
  \varphi_i &=& \mu(q^i-1)(1-q^{d-i+1})(\mu^* q^{i-1} - h^* q^{d-i})
        \qquad  (1 \leq i \leq d).
     \label{eq:Iaux2}
\end{eqnarray}
Evaluating the equation on the left in (\ref{eq:ai}) using
(\ref{eq:Iths}), (\ref{eq:Iaux1}), (\ref{eq:Iaux2}) 
we routinely find $a_i=\eta$ for $0 \leq i \leq d$.

(i)$\Rightarrow$(iii):
Setting $h+\mu=0$ in Lemma \ref{lem:Igammagammas} we find
$\theta_i+\theta_{d-i}=2 \eta$ for $0 \leq i \leq d$.
Setting $\tau=0$ and $h+\mu=0$ in Lemma \ref{lem:Ivphiphi}
we find $\varphi_i=- \phi_i$ for $1 \leq i \leq d$.

(ii)$\Rightarrow$(i): We have $\tau=0$ and $(h+\mu)(h^*+\mu^*)=0$ by
Theorem \ref{thm:I}. Suppose $h+\mu \neq 0$. Then we must have
$h^*+\mu^*=0$, so that Lemma \ref{lem:Ia1a0} implies $\mu^* (h+\mu)=0$.
Observe that we have $\mu^* \neq 0$; otherwise $h^*=h^*+\mu^*=0$ so that
$\theta^*_0=\theta^*_1$. Hence $h+\mu=0$.

(iii)$\Rightarrow$(i): 
Consider the quantity $\theta_0+\theta_d - \theta_1 - \theta_{d-1}$.
By assumption this quantity is $0$.
By Lemma \ref{lem:Igammagammas} this quantity is
$(q-1)(q^{d-1}-1)(h+\mu)$ so $h+\mu=0$.
Setting $\varphi_i+\phi_i=0$, $h+\mu=0$ in Lemma \ref{lem:Ivphiphi}
we find $2\tau=0$.
Observe $\text{\rm Char}(\mathbb{K}) \neq 2$; otherwise $\theta_d=\theta_0$
by Lemma \ref{lem:Igammagammas}. We conclude $\tau=0$.
\end{proof}

\medskip

\begin{theorem}        \label{thm:Isaux}
Assume $d \geq 3$, $q \neq 1$, $q \neq -1$.
Then the following are equivalent.
\begin{itemize}
\item[(i)] $\tau=0$ and $h^*+\mu^*=0$.
\item[(ii)] $a^*_i$ is independent of $i$ for $0 \leq i \leq d$.
\item[(iii)] $\theta^*_i + \theta^*_{d-i}$ is independent of $i$
   for $0 \leq i \leq d$, and
   $\varphi_i = -\phi_{d-i+1}$ for $1 \leq i \leq d$.
\end{itemize}
Suppose (i)--(iii) hold. Then the common value of 
$a^*_i$ is $\eta^*$, and the common value of
$\theta^*_i + \theta^*_{d-i}$ is $2 \eta^*$.
\end{theorem}

\begin{proof}
Similar to the proof of Thoerem \ref{thm:Iaux}.
\end{proof}

\section{Case II: $d \geq 3$, $q=1$, $\text{\rm Char}(\mathbb{K}) \neq 2$}

In this section we assume $d \geq 3$, $q=1$, $\text{\rm Char}(\mathbb{K}) \neq 2$.

\medskip

\begin{theorem}   \cite{T:array}
There exist scalars
$\eta$, $\mu$, $h$, $\eta^*$, $\mu^*$, $h^*$, $\tau$ in 
$\overline{\mathbb{K}}$ such that for $0 \leq i \leq d$
\begin{eqnarray}
\theta_i &=& \eta + \mu(i-d/2)+hi(d-i),
   \label{eq:IIth}  \\
\theta^*_i &=& \eta^* + \mu^*(i-d/2)+h^* i(d-i),
   \label{eq:IIths} 
\end{eqnarray}
and for $1 \leq i \leq d$
\begin{eqnarray}
\varphi_i &=& i(d-i+1)(\tau-\mu\mu^*/2+(h\mu^*+\mu h^*)(i-(d+1)/2)+hh^*(i-1)(d-i)),
   \label{eq:IIphi}  \\
\phi_i &=& i(d-i+1)(\tau+\mu\mu^*/2+(h\mu^*-\mu h^*)(i-(d+1)/2)+hh^*(i-1)(d-i)).
   \label{eq:IIphi2}
\end{eqnarray}
\end{theorem}

\begin{proof}
These are  (35), (36), (38), (39) in \cite{T:array}
after a change of variables.
\end{proof}

\medskip

\begin{remark}          \label{rem:II}
If $h=0$ then $\mu \neq 0$, otherewise $\theta_{1}=\theta_{0}$.
Similarly if $h^*=0$ then $\mu^* \neq 0$. 
For any prime $i$ such that $i \leq d$ we have 
$\text{\rm Char}(\mathbb{K}) \neq i$; otherwise $\varphi_i=0$ by
(\ref{eq:IIphi}).
\end{remark}

\medskip

\begin{lemma}         \label{lem:IIada0}
$H = 2 \tau +hh^*(d-1)^2$.
\end{lemma}

\begin{proof}
It is routine to verify this equation
using (\ref{eq:ai}), (\ref{eq:IIth}), (\ref{eq:IIths}) and
(\ref{eq:IIphi}).
\end{proof}

\medskip

\begin{lemma}         \label{lem:IItau}
Assume $H=0$. Then
\begin{equation}       \label{eq:IItau}
\tau = -hh^*(d-1)^2/2.
\end{equation}
\end{lemma}

\begin{proof}
Follows from Lemma \ref{lem:IIada0}.
\end{proof}

\medskip

\begin{lemma}               \label{lem:IIadm1a1}
Assume $H=0$. Then the following coincide.
\[
    \frac{(a_1-a_{d-1})(\theta^*_0-\theta^*_3)(\theta^*_{d-3}-\theta^*_d)}
         {\theta^*_0 - \theta^*_d},
\]
\[
    \frac{(a^*_1-a^*_{d-1})(\theta_0-\theta_3)(\theta_{d-3}-\theta_d)}
         {\theta_0 - \theta_d},
\]
\[
   -36d^{-1}(d-1)(d-2)hh^*.
\]
\end{lemma}

\begin{proof}
It is routine to verify the coincidence
using (\ref{eq:ai}), (\ref{eq:IIth}), (\ref{eq:IIths}), (\ref{eq:IIphi})
and (\ref{eq:IItau}).
\end{proof}

\medskip

\begin{lemma}             \label{lem:IIequiv1}
Assume $H=0$. Then the following are equivalent.
\begin{itemize}
\item[(i)] $a_1=a_{d-1}$,
\item[(ii)] $a^*_1 = a^*_{d-1}$,
\item[(iii)] $hh^*=0$,
\item[(iv)] $\tau=0$.
\end{itemize}
\end{lemma}

\begin{proof}
Follows from Lemma \ref{lem:IIadm1a1} and Remark \ref{rem:II}.
\end{proof}

\medskip

\begin{theorem}        \label{thm:II}
Assume $d \geq 3$, $q=1$, $\text{\rm Char}(\mathbb{K}) \neq 2$.
Then  the following are equivalent.
\begin{itemize}
\item[(i)] $a_0=a_d$ and $a_1=a_{d-1}$,
\item[(ii)] $a^*_0=a^*_d$ and $a^*_1=a^*_{d-1}$,
\item[(iii)] $a_i=a_{d-i}$ and $a^*_i=a^*_{d-i}$ for $0 \leq i \leq d$,
\item[(iv)] $hh^*=0$ and $\tau=0$.
\end{itemize}
\end{theorem}

\begin{proof}
The conditions (i), (ii), (iv) are equivalent by Lemmas
\ref{lem:IIada0} and \ref{lem:IIequiv1}. 
Clearly (iii) implies (i).  
We show (iv) implies (iii).
It is routine to verify $a_i-a_{d-i}=0$ and $a^*_{i}-a^*_{d-i}=0$
for each case of $h=0$, $h^*=0$ by using
(\ref{eq:ai}), (\ref{eq:IIth}), (\ref{eq:IIths}), (\ref{eq:IIphi})
with $\tau=0$. 
\end{proof}

\medskip

\begin{lemma}          \label{lem:IIgammagammas}
For $0 \leq i \leq d$,
\begin{eqnarray*}
  \theta_i +\theta_{d-i} &=& 2(\eta+hi(d-i)), \\
  \theta^*_i+\theta^*_{d-i} &=& 2(\eta^*+h^* i(d-i)). 
\end{eqnarray*}
\end{lemma}

\begin{proof}
It is routine to verify these equations using
(\ref{eq:IIth}) and (\ref{eq:IIths}).
\end{proof}

\medskip

\begin{lemma}              \label{lem:IIvphiphi}
For $1 \leq i \leq d$,
\begin{eqnarray*}
  \varphi_i + \phi_i 
  &=& i(d-i+1)(2\tau-(d-2i+1)h\mu^*+2hh^*(d-i)(i-1)), \\
  \varphi_i + \phi_{d-i+1}
 &=&  i(d-i+1)(2\tau-(d-2i+1)h^*\mu+2hh^*(d-i)(i-1)).
\end{eqnarray*}
\end{lemma}

\begin{proof}
It is routine to verify these equations using
(\ref{eq:IIphi}) and (\ref{eq:IIphi2}).
\end{proof}

\medskip

\begin{lemma}          \label{lem:IIa1a0}
The following hold.
\begin{itemize}
\item[(i)] Assume $\tau=0$ and $h^* = 0$. Then
\[
a_0-a_1 = 2(d-1)h.
\]
\item[(ii)] Assume $\tau=0$ and $h = 0$. Then
\[
a^*_0-a^*_1 = 2(d-1) h^*.
\]
\end{itemize}
\end{lemma}

\begin{proof}
It is routine to verify these equations using
(\ref{eq:ai}), (\ref{eq:IIth}), (\ref{eq:IIths}), (\ref{eq:IIphi}).
\end{proof}

\medskip

\begin{theorem}        \label{thm:IIaux}
Assume $d \geq 3$, $q=1$, $\text{\rm Char}(\mathbb{K}) \neq 2$.
Then the following are equivalent.
\begin{itemize}
\item[(i)] $h=0$ and $\tau=0$.
\item[(ii)] $a_i$ is independent of $i$ for $0 \leq i \leq d$.
\item[(iii)] $\theta_i + \theta_{d-i}$ is independent of $i$
   for $0 \leq i \leq d$, and
   $\varphi_i = - \phi_i$ for $1 \leq i \leq d$.
\end{itemize}
Suppose (i)--(iii) hold. Then the common value of 
$a_i$ is $\eta$, and the common value of
$\theta_i + \theta_{d-i}$ is $2 \eta$.
\end{theorem}

\begin{proof}
(i)$\Rightarrow$(ii):
Evaluating (\ref{eq:IIth}), (\ref{eq:IIphi})
using $h=0$ and $\tau=0$ we find
\begin{eqnarray}
 \theta_i &=& \eta + (i-d/2)\mu     \qquad (0 \leq i \leq d),
    \label{eq:IIaux1}  \\
  \varphi_i &=& - i(d-i+1)\mu(\mu^*+h^*(d-2i+1))/2
      \qquad (1 \leq i \leq d).
    \label{eq:IIaux2}
\end{eqnarray}
Evaluating the equation on the left in (\ref{eq:ai}) using
(\ref{eq:IIths}), (\ref{eq:IIaux1}), (\ref{eq:IIaux2}) we routinely
find $a_i=\eta$ for $0 \leq i \leq d$.

(i)$\Rightarrow$(iii):
Setting $h=0$ in Lemma \ref{lem:IIgammagammas} we find
$\theta_i+\theta_{d-i} = 2 \eta$ for $0 \leq i \leq d$.
Setting $h=0$ and $\tau=0$ in Lemma \ref{lem:IIvphiphi}
we find $\varphi_i = - \phi_i$ for $1 \leq i \leq d$.

(ii)$\Rightarrow$(i): We have $\tau=0$ and $hh^*=0$ by Theorem \ref{thm:II}.
Suppose $h \neq 0$. Then we must have $h^*=0$. Then Lemma \ref{lem:IIa1a0}
implies $h=0$.

(iii)$\Rightarrow$(i):
Consider the quantity $\theta_0+\theta_d-\theta_1-\theta_{d-1}$.
By assumption this quantity is $0$. By Lemma \ref{lem:IIgammagammas}
this quantity is $2(1-d)h$ so $h=0$.
Setting $\varphi_i+\phi_i=0$, $h=0$ in Lemma \ref{lem:IIvphiphi}
we find $\tau=0$.
\end{proof}

\medskip

\begin{theorem}        \label{thm:IIsaux}
Assume $d \geq 3$, $q=1$, $\text{\rm Char}(\mathbb{K}) \neq 2$.
Then the following are equivalent.
\begin{itemize}
\item[(i)] $h^*=0$ and $\tau=0$.
\item[(ii)] $a^*_i$ is independent of $i$ for $0 \leq i \leq d$.
\item[(iii)] $\theta^*_i + \theta^*_{d-i}$ is independent of $i$
   for $0 \leq i \leq d$, and
   $\varphi_i = - \phi_{d-i+1}$ for $1 \leq i \leq d$.
\end{itemize}
Suppose (i)--(iii) hold. Then the common value of 
$a^*_i$ is $\eta^*$, and the common value of
$\theta^*_i + \theta^*_{d-i}$ is $2 \eta^*$.
\end{theorem}

\begin{proof}
Similar to the proof of Theorem \ref{thm:IIaux}.
\end{proof}

\section{Case III: $d \geq 3$, $q=-1$, $\text{\rm Char}(\mathbb{K}) \neq 2$, $d$ even}

In this section we assume $d \geq 3$,  $q=-1$, $\text{\rm Char}(\mathbb{K}) \neq 2$,
and $d$ is even.

\medskip

\begin{theorem}     \cite[Theorem 5.16, Example 5.14]{T:array}
There exist scalars $\eta$, $h$, $s$, $\eta^*$, $h^*$, $s^*$,
$\tau$ in $\overline{\mathbb{K}}$ such that for $0 \leq i \leq d$
\begin{eqnarray}
\theta_i &=&
  \begin{cases}
     \eta+s+h(i-d/2)  & \text{\rm if $i$ is even}, \\
     \eta-s -h(i-d/2) & \text{\rm if $i$ is odd},
  \end{cases}   
         \label{eq:IIIeventh}   \\
\theta^*_i &=&
   \begin{cases}
     \eta^* +s^* +h^*(i-d/2)  &   \text{\rm if $i$ is even}, \\
     \eta^* - s^* -h^*(i-d/2) &  \text{\rm if $i$ is odd},
   \end{cases}
      \label{eq:IIIevenths}
\end{eqnarray}
and for $1 \leq i \leq d$
\begin{eqnarray}    
\varphi_i &=&
   \begin{cases}
      i(\tau-sh^*-s^*h-hh^*(i-(d+1)/2)) &   \text{\rm if $i$ is even}, \\
      (d-i+1)(\tau+sh^*+s^*h+hh^*(i-(d+1)/2))  & \text{\rm if $i$ is odd},
   \end{cases}
     \label{eq:IIIevenphi}  \\
\phi_i &=&
   \begin{cases}
      i(\tau-sh^*+s^*h+hh^*(i-(d+1)/2)) &   \text{\rm if $i$ is even}, \\
      (d-i+1)(\tau+sh^*-s^*h - hh^*(i-(d+1)/2))  & \text{\rm if $i$ is odd}.
   \end{cases}
    \label{eq:IIIevenphi2}
\end{eqnarray}
\end{theorem}

\begin{proof}
These are  (19)--(22) in \cite{T:array}
after a change of variables.
\end{proof}

\medskip

\begin{remark}        \label{rem:IIIeven}
We have $h \neq 0$; otherwise $\theta_0=\theta_2$ by (\ref{eq:IIIeventh}).
Similary we have $h^* \neq 0$.
For any prime $i$ such that $i \leq d/2$ we have 
$\text{\rm Char}(\mathbb{K})\neq i$; otherwise $\varphi_{2i}=0$
by ($\ref{eq:IIIevenphi}$). By this and since $\text{\rm Char}(\mathbb{K})\neq 2$
we find $\text{\rm Char}(\mathbb{K})$ is either $0$ or an odd prime
greater than $d/2$.
Observe neither of $d$, $d-2$ vanish in $\mathbb{K}$ since otherwise
$\text{\rm Char}(\mathbb{K})$ must divide $d/2$ or $(d-2)/2$.
\end{remark}

\medskip

\begin{lemma}     \label{lem:IIIevenada0}
$H = 2(d-1)\tau+4ss^*$.
\end{lemma}

\begin{proof}
It is routine to verify this equation
using (\ref{eq:ai}), (\ref{eq:IIIeventh}), (\ref{eq:IIIevenths}),
(\ref{eq:IIIevenphi}).
\end{proof}

\medskip

\begin{lemma}       \label{lem:IIIevenr1}
Assume $H=0$. Then $d-1$ is nonzero in $\mathbb{K}$ and
\begin{equation}    \label{eq:IIIevenr1}
  \tau = \frac{2ss^*}{1-d}.
\end{equation}
\end{lemma}

\begin{proof}
Suppose $d-1$ is zero in $\mathbb{K}$. Then Lemma \ref{lem:IIIevenada0} implies
$ss^*=0$. If $s=0$ then
$\theta_1=\theta_0$ by (\ref{eq:IIIeventh}).
If $s^*=0$ then $\theta^*_1=\theta^*_0$ by (\ref{eq:IIIevenths}).
Hence $d-1$ is nonzero and (\ref{eq:IIIevenr1}) follows.
\end{proof}

\medskip

\begin{lemma}          \label{lem:IIIevenadm1a1}
Assume $H=0$. Then the following coincide.
\[
  \frac{(a_1-a_{d-1})(\theta^*_0 - \theta^*_3)(\theta^*_{d-3} - \theta^*_{d})}
       {\theta^*_0-\theta^*_d},
\]
\[
  \frac{(a^*_1-a^*_{d-1})(\theta_0 - \theta_3)(\theta_{d-3} - \theta_{d})}
       {\theta_0-\theta_d},
\]
\[
   \frac{16(d-2)ss^*}{d(d-1)}.
\]
\end{lemma}

\begin{proof} 
It is routine to verify the coincidence using (\ref{eq:ai}),
(\ref{eq:IIIeventh}), (\ref{eq:IIIevenths}), (\ref{eq:IIIevenphi})
and (\ref{eq:IIIevenr1}).
\end{proof}

\medskip

\begin{lemma}    \label{lem:IIIevenequiv1}
Assume $H=0$. Then the following are equivalent.
\begin{itemize}
\item[(i)] $a_1=a_{d-1}$,
\item[(ii)] $a^*_1 = a^*_{d-1}$,
\item[(iii)] $ss^*=0$,
\item[(iv)] $\tau=0$.
\end{itemize}
\end{lemma}

\begin{proof}
Follows from Lemma \ref{lem:IIIevenadm1a1} and (\ref{eq:IIIevenr1}).
\end{proof}

\medskip

\begin{theorem}    \label{thm:IIIeven}
\samepage
Assume $d \geq 3$, $q=-1$, $\text{\rm Char}(\mathbb{K}) \neq 2$, 
and $d$ is even.
Then the following are equivalent.
\begin{itemize}
\item[(i)]  $a_0 = a_d$ and $a_1=a_{d-1}$,
\item[(ii)] $a^*_0=a^*_d$ and $a^*_1=a^*_{d-1}$,
\item[(iii)] $a_i=a_{d-i}$ and $a^*_i=a^*_{d-i}$ for $0 \leq i\leq d$,
\item[(iv)] $ss^*=0$ and $\tau=0$.
\end{itemize}
\end{theorem}

\begin{proof}
The conditions (i), (ii), (iv) are equivalent by Lemmas
\ref{lem:IIIevenada0} and \ref{lem:IIIevenequiv1}. 
Clearly (iii) implies (i).
We show (iv) implies (iii).
It is routine to verify $a_{i}-a_{d-i}=0$ and
$a^*_{i}- a^*_{d-i}=0$ for each case of $s=0$ and $s^*=0$
using (\ref{eq:ai}), (\ref{eq:IIIeventh}),
(\ref{eq:IIIevenths}), (\ref{eq:IIIevenphi}) with
$\tau=0$.
\end{proof}

\medskip

\begin{lemma}          \label{lem:IIIevengamma}
For $0 \leq i \leq d$,
\begin{eqnarray*}
   \theta_i+\theta_{d-i} 
       &=& \begin{cases}
             2(\eta+s) & \text{\rm if $i$ is even}, \\
             2(\eta-s) & \text{\rm if $i$ is odd},
           \end{cases}  \\
   \theta^*_i+\theta^*_{d-i} 
       &=& \begin{cases}
            2(\eta^*+s^*) &  \text{\rm if $i$ is even}, \\
            2(\eta^*-s^*) & \text{\rm if $i$ is odd}.
           \end{cases}
\end{eqnarray*}
\end{lemma}

\begin{proof}
It is routine to verify these equations using
(\ref{eq:IIIeventh}) and (\ref{eq:IIIevenths}).
\end{proof}

\medskip

\begin{lemma}             \label{lem:IIIevenphi}
For $1 \leq i \leq d$,
\begin{eqnarray*}
\varphi_i + \phi_i 
 &=&  \begin{cases}
        2i(\tau-sh^*)  &  \text{\rm if $i$ is even}, \\
        2(d-i+1)(\tau+sh^*) & \text{\rm if $i$ is odd}, 
      \end{cases}   \\
\varphi_i + \phi_{d-i+1}
 &=&  \begin{cases}
        2i(\tau-s^*h)  &  \text{\rm if $i$ is even}, \\
        2(d-i+1)(\tau+s^*h) & \text{\rm if $i$ is odd}.
      \end{cases}  
\end{eqnarray*}
\end{lemma}

\begin{proof}
It is routine to verify these equations using
(\ref{eq:IIIevenphi}) and (\ref{eq:IIIevenphi2}).
\end{proof}

\medskip

\begin{lemma}       \label{lem:IIIevena0a1}
The following hold.
\begin{itemize}
\item[(i)] Assume $\tau=0$ and $s^*=0$.  Then
each of $d-1$, $d-3$ is nonzero in $\mathbb{K}$ and
\[
   a_0 - a_1 = 
    \frac{4s}{(d-1)(d-3)}.
\]
\item[(ii)] Assume $\tau=0$ and $s=0$. Then
each of $d-1$, $d-3$ is nonzero in $\mathbb{K}$ and
\[
  a^*_0 - a^*_1 = 
    \frac{4s^*}{(d-1)(d-3)}.
\]

\end{itemize}
\end{lemma}

\begin{proof}
We show (i). 
If $d-1$ is zero in $\mathbb{K}$ then $\theta^*_0=\theta^*_1$ by (\ref{eq:IIIevenths}).
If $d-3$ is zero in $\mathbb{K}$ then $\theta^*_0=\theta^*_3$ by (\ref{eq:IIIevenths}).
Hence each of $d-1$, $d-3$ is nonzero.
Now we routinely find the equation for $a_0-a_1$
using (\ref{eq:ai}), (\ref{eq:IIIeventh}), (\ref{eq:IIIevenths})
and (\ref{eq:IIIevenphi}).
The proof of (ii) is similar.
\end{proof}

\medskip

\begin{theorem}        \label{thm:IIIaux}
Assume $d \geq 3$, $q=-1$, $\text{\rm Char}(\mathbb{K}) \neq 2$, 
and $d$ is even.
Then the following are equivalent.
\begin{itemize}
\item[(i)] $s=0$ and $\tau=0$.
\item[(ii)] $a_i$ is independent of $i$ for $0 \leq i \leq d$.
\item[(iii)] $\theta_i + \theta_{d-i}$ is independent of $i$
   for $0 \leq i \leq d$, and
   $\varphi_i = -\phi_i$ for $1 \leq i \leq d$.
\end{itemize}
Suppose (i)--(iii) hold. Then the common value of 
$a_i$ is $\eta$, and the common value of
$\theta_i + \theta_{d-i}$ is $2\eta$.
\end{theorem}

\begin{proof}
(i)$\Rightarrow$(ii):
Evaluating (\ref{eq:IIIeventh}), 
(\ref{eq:IIIevenphi}) using $s=0$ and $\tau=0$ we find
\begin{eqnarray}
 \theta_i &=&
   \begin{cases}
      \eta+h(i-d/2) & \text{ if $i$ is even},  \\
      \eta-h(i-d/2) & \text{ if $i$ is odd},
   \end{cases}
     \label{eq:IIIaux1}  \\
  \varphi_i &=&
   \begin{cases}
      -hi(s^*+h^*(i-(d+1)/2)) &  \text{ if $i$ is even},  \\
      h(d-i+1)(s^*+h^*(i-(d+1)/2))  & \text{ if $i$ is odd}.
   \end{cases}
      \label{eq:IIIaux3}
\end{eqnarray}
Evaluating the equation on the left in (\ref{eq:ai}) using
(\ref{eq:IIIaux1}),  (\ref{eq:IIIaux3}) we routinely
find $a_i=\eta$ for $0 \leq i \leq d$.

(i)$\Rightarrow$(iii):
Setting $s=0$ in Lemma \ref{lem:IIIevengamma}
we find $\theta_i+\theta_{d-i}=2\eta$ for $0 \leq i \leq d$.
Setting $s=0$ and $\tau=0$ in Lemma \ref{lem:IIIevenphi}
we find $\varphi_i = - \phi_i$ for $1 \leq i \leq d$.

(ii)$\Rightarrow$(i): Suppose (i) does not hold. Then from
Theorem \ref{thm:IIIeven}, we must have
$s^*=0$ and $\tau=0$.
From our assuption, we have $a_1-a_0=0$, so Lemma \ref{lem:IIIevena0a1}
implies $s=0$, a contradiction.

(iii)$\Rightarrow$(i):
From Lemma \ref{lem:IIIevenphi} for $i=1,2$,
\[
  0 = \varphi_1 + \phi_1 = 2d(\tau+sh^*),
\]
\[
  0 = \varphi_2 + \phi_2 = 4(\tau-sh^*).
\]
These equations imply $\tau=0$ and $s=0$.
\end{proof}

\medskip

\begin{theorem}        \label{thm:IIIsaux}
Assume $d \geq 3$, $q=-1$, $\text{\rm Char}(\mathbb{K}) \neq 2$, 
and $d$ is even.
Then the following are equivalent.
\begin{itemize}
\item[(i)]  $s^*=0$ and $\tau=0$.
\item[(ii)] $a^*_i$ is independent of $i$ for $0 \leq i \leq d$.
\item[(iii)] $\theta^*_i + \theta^*_{d-i}$ is independent of $i$
   for $0 \leq i \leq d$, and
   $\varphi_i = - \phi_{d-i+1}$ for $1 \leq i \leq d$.
\end{itemize}
Suppose (i)--(iii) hold. Then the common value of 
$a^*_i$ is $\eta^*$, and the common value of
$\theta^*_i + \theta^*_{d-i}$ is $2\eta^*$.
\end{theorem}

\begin{proof}
Similar to the proof of Theorem \ref{thm:IIIaux}.
\end{proof}

\section{Case IV: $d \geq 3$, $q=-1$, $\text{\rm Char}(\mathbb{K}) \neq 2$, $d$ odd}

In this section we assume $d \geq 3$, $q=-1$, $\text{\rm Char}(\mathbb{K}) \neq 2$,
and $d$ is odd.

\medskip

\begin{theorem}  \cite[Theorem 5.16, Example 5.14]{T:array}
There exist scalars $\eta$, $h$, $s$, $\eta^*$, $h^*$, $s^*$, $\tau$
in $\overline{\mathbb{K}}$ such that for $0 \leq i \leq d$
\begin{eqnarray}
\theta_i &=&
  \begin{cases}
     \eta+s+h(i-d/2)  & \text{\rm if $i$ is even}, \\
     \eta-s-h(i-d/2) & \text{\rm if $i$ is odd},
  \end{cases}   
         \label{eq:IIIoddth}   \\
\theta^*_i &=&
   \begin{cases}
     \eta^* + s^* +h^*(i-d/2)  &   \text{\rm if $i$ is even}, \\
     \eta^*-s^* -h^*(i-d/2) &  \text{\rm if $i$ is odd},
   \end{cases}
      \label{eq:IIIoddths}
\end{eqnarray}
and for $1 \leq i \leq d$
\begin{eqnarray}    
\varphi_i &=&
   \begin{cases}
      hh^* i(d-i+1) &   \text{\rm if $i$ is even}, \\
      \tau-2ss^*+i(d-i+1)hh^*-2(hs^*+h^*s)(i-(d+1)/2)  & \text{\rm if $i$ is odd},
   \end{cases}
      \label{eq:IIIoddphi}  \\
\phi_i &=&
   \begin{cases}
      hh^* i(d-i+1) &   \text{\rm if $i$ is even}, \\
      \tau+2ss^*+i(d-i+1)hh^*-2(hs^*-h^*s)(i-(d+1)/2)  & \text{\rm if $i$ is odd}.
   \end{cases}
      \label{eq:IIIoddphi2}
\end{eqnarray}
\end{theorem}

\medskip

\begin{remark}         \label{rem:IIIodd}
Observe $hh^* \neq 0$, and $\text{\rm Char}(\mathbb{K})$ is
either $0$ or an odd prime greater than $d/2$.
Also observe $d-1$ does not vanish in $\mathbb{K}$.
These can be observed in a similar way as Remark \ref{rem:IIIeven}.
\end{remark}

\medskip

\medskip

\begin{lemma}     \label{lem:IIIoddada0}
$H= 2\tau+(d^2+1)hh^*$.
\end{lemma}

\begin{proof}
It is routine to verify this equation using (\ref{eq:ai}), (\ref{eq:IIIoddth}),
(\ref{eq:IIIoddths}) and (\ref{eq:IIIoddphi}).
\end{proof}

\medskip

\begin{lemma}       \label{lem:IIIoddnonzero}
Assume $H=0$. Then
\begin{equation}        \label{eq:IIIoddss}
   \tau=-(d^2+1)hh^*/2.
\end{equation}
\end{lemma}

\begin{proof}
Follows from Lemma \ref{lem:IIIoddada0}.
\end{proof}

\medskip

\begin{lemma}         \label{lem:IIIoddadm1a1}
Assume $H=0$. Then the following coincide.
\[
   \frac{(a_1-a_{d-1})(\theta^*_0-\theta^*_3)(\theta^*_{d-3}-\theta^*_d)}
        {\theta^*_0-\theta^*_d},
\]
\[
   \frac{(a^*_1-a^*_{d-1})(\theta_0-\theta_3)(\theta_{d-3}-\theta_d)}
        {\theta_0-\theta_d},
\]
\[
   -4(d-1) h h^*.
\]
\end{lemma}

\begin{proof}
It is routine to verify the coincidence using
(\ref{eq:ai}), (\ref{eq:IIIoddth}), (\ref{eq:IIIoddths}),
(\ref{eq:IIIoddphi}) and (\ref{eq:IIIoddss}).
\end{proof}

\medskip

\begin{theorem}        \label{thm:IIIodd}
Assume $d \geq 3$, $q=-1$, $\text{\rm Char}(\mathbb{K}) \neq 2$, and
$d$ is odd.
If $a_0=a_d$ then $a_1 \neq a_{d-1}$ and $a^*_1 \neq a^*_{d-1}$.
\end{theorem}

\begin{proof}
Follows from Lemma \ref{lem:IIIoddadm1a1} and Remark \ref{rem:IIIodd}.
\end{proof}

\medskip

\begin{theorem}         \label{thm:IIIoddaux}
Assume $d \geq 3$, $q=-1$, $\text{\rm Char}(\mathbb{K}) \neq 2$, and
$d$ is odd. Then $\varphi_2+\phi_2 \neq 0$ and
$\varphi_2 + \phi_{d-1} \neq 0$.
\end{theorem}

\begin{proof}
From (\ref{eq:IIIoddphi}) and (\ref{eq:IIIoddphi2}),
$\varphi_2 + \phi_2 = \varphi_2 + \phi_{d-1} = 4(d-1)hh^* \neq 0$
by Remark \ref{rem:IIIodd}.
\end{proof}

\section{Case V: $d \geq 3$, $q=1$, $\text{\rm Char}(\mathbb{K}) = 2$}

In this section we assume $d \geq 3$, $q=1$, $\text{\rm Char}(\mathbb{K}) = 2$.

\medskip

\begin{theorem}  \cite[Theorem 5.16, Example 5.15]{T:array}
We have $d=3$, and there exist scalars
$h$, $s$, $h^*$, $s^*$, $r$ in $\overline{\mathbb{K}}$ such that
\[
 \begin{array}{lll}
    \theta_1 = \theta_0 + h(s+1),
  & \theta_2 = \theta_0 + h,
  & \theta_3 = \theta_0 + hs,
  \\
    \theta^*_1 = \theta^*_0 + h^*(s^*+1), 
  & \theta^*_2 = \theta^*_0 + h^*,
  & \theta^*_3 = \theta^*_0 + h^*s^*,
  \\
    \varphi_1 = hh^* r,
  & \varphi_2 = hh^*, 
  & \varphi_3 = hh^*(r+s+s^*), 
  \\
    \phi_1 = hh^*(r+s(1+s^*)),
  & \phi_2 = hh^*,
  & \phi_3 = hh^*(r+s^*(1+s)).
  \end{array}
\]
\end{theorem}

\medskip

\begin{remark}
Each of $h$, $h^*$, $s$, $s^*$ is nonzero, and each of
$s$, $s^*$ is not equal to $1$.
\end{remark}

\medskip

\begin{lemma}        \label{lem:IVada0}
\begin{eqnarray}    
   a_0-a_3 &=& \frac{h s^*(1+s)}{1+s^*},
       \label{eq:IVada0} \\
   a^*_0-a^*_3 &=& \frac{h^* s(1+s^*)}{1+s}.
       \label{eq:IVasdas0} 
\end{eqnarray}
\end{lemma}

\begin{proof}
Obtained by a routine computation. We remark that $2=0$ and $1=-1$ since
$\text{\rm Char}(\mathbb{K})=2$.
\end{proof}

\medskip

\begin{theorem}  \label{thm:IV}
Assume $d \geq 3$, $q=1$, $\text{\rm Char}(\mathbb{K})=2$.
Then $a_0 \neq a_d$ and $a^*_0 \neq a^*_d$.
\end{theorem}

\begin{proof}
Immediate from Lemma \ref{lem:IVada0} and since none of 
$h$, $h^*$, $s$, $s^*$, $1+s$, $1+s^*$ is zero.
\end{proof}

\medskip

\begin{theorem}   \label{thm:IVaux}
Assume $d \geq 3$, $q=1$, $\text{\rm Char}(\mathbb{K})=2$.
Then $\varphi_1+\phi_1 \neq 0$ and $\varphi_1+\phi_d \neq 0$.
\end{theorem}

\begin{proof}
We have
\begin{eqnarray*}
 \varphi_1 + \phi_1 &=& hh^*s(1+s^*), \\
 \varphi_1 + \phi_3 &=& hh^* s^*(1+s).
\end{eqnarray*}
These values are nonzero since none of $h$, $h^*$, $s$, $s^*$,
$1+s$, $1+s^*$ is zero.
\end{proof}

\bigskip

\bibliographystyle{plain}

\bigskip\bigskip\noindent
Kazumasa Nomura\\
College of Liberal Arts and Sciences\\
Tokyo Medical and Dental University\\
Kohnodai, Ichikawa, 272-0827 Japan\\
email: nomura.las@tmd.ac.jp

\bigskip\noindent
Paul Terwilliger\\
Department of Mathematics\\
University of Wisconsin\\
480 Lincoln drive, Madison, Wisconsin, 53706 USA\\
email: terwilli@math.wisc.edu

\bigskip\noindent
{\bf Keywords.}
Leonard pair, Terwilliger algebra, Askey scheme,
$q$-Racah polynomial.

\noindent
{\bf 2000 Mathematics Subject Classification}.
05E30, 05E35, 33C45, 33D45.

\end{document}